\documentclass[11pt]{article}
\usepackage[paperheight=30cm,paperwidth=20cm,textwidth=15cm]{geometry}
\usepackage{amsfonts}
\usepackage{amsmath}
\usepackage{mathrsfs}
\usepackage{caption}
\usepackage{multicol}
\usepackage{blindtext}
\usepackage{multicol}
\usepackage{amssymb}
\usepackage{graphicx}
\usepackage{pdfpages}
\def\qed{\hfill $\square$ \\}

\def\Aut{{\rm Aut}}

\newtheorem{theorem}{Theorem}[section]
\newtheorem{lemma}[theorem]{Lemma}%
\newtheorem{corollary}[theorem]{Corollary}%
\newtheorem{question}[theorem]{Question}%
\def\pf{\noindent{\it Proof.} }
 
\begin{document}



\title{{\bf On automorphism groups of binary cyclic codes}}


\author{Jicheng Ma \\[+3pt]
Chongqing Key Lab. of Group {\&} Graph Theories and Applications,\\ Chongqing University of Arts and Sciences,\\
Chongqing 402160, China \\[+3pt] {\normalsize ma\_jicheng@hotmail.com}\\[+10pt]
Guiying Yan \\[+3pt]
Academy of Mathematics and System Science,\\ Chinese Academy of Sciences,\\  Beijing 100080, China \\[+3pt] {\normalsize yangy@amss.ac.cn}\\[+10pt] \\}

\date{}

\maketitle

{\bf Abstract}:
A lot of attention has been paid to the investigation of the algebraic properties of linear codes. In most cases, this investigation involves the determination of required code automorphisms, which are useful for decoders, such as the automorphism ensemble decoder. 
It is worth noting that the examination of the automorphism groups of discrete symmetric objects has long been a highly regarded field of research.
Cyclic codes, as a significant subclass of linear codes, can be constructed and analyzed using algebraic methods. 
And due to its cyclic nature, they have efficient encoding and decoding algorithms. 
To date, cyclic codes have found applications in various domains, including consumer electronics, data storage systems, and communication systems. 
In this paper, we investigate the full automorphism groups of binary cyclic codes. In particular, we present constructions for binary cyclic codes of long lengths that facilitate the determination of the full automorphism groups.
 
 \medskip



 \section{Introduction}
 
   In coding theory, establishing a natural connection between codes and groups is of great interest, as  a group acting on a code may offer valuable insights into the structure of the code. This concept has been widely recognized. 
The Mathieu groups acting on the Golay codes \cite{Conway}, the affine linear group ${\rm AGL(n, 2)}$ acting on the Reed-Muller codes \cite{MacWilliams}, and the general linear group ${\rm GL}(n, 2)$ acting on Hamming codes are examples that demonstrate this connection. 
In particular, determining the automorphism group of a linear code is important due to both theoretical interest in itself and applications in decoding algorithms. 
So far, 
the complete determination of the automorphism group of only few classes of linear codes are classified, such as Hamming codes, Reed-Muller codes, etc. For more detailed information, we refer to the Handbook article by Huffman \cite[Ch. 17]{Pless}.

In 1987, D{\"u}r \cite{Dur} investigated the automorphism group of Reed-Solomon and extended Reed-Solomon codes. 
In 1993, Berger and Charpin \cite{BergerCharpin} proved that the automorphism group of generalized Reed-Muller codes is the general linear nonhomogeneous group, and in \cite{BergerCharpin1999} they also characterized that the automorphism group of Bose-Chaudhuri-Hocquenghem (BCH) codes and some affine-invariant codes.   
Subsequently, in 1996,  Berger \cite{Berger} investigated the automorphism groups of non-trivial affine-invariant codes and established that these automorphism groups are always subgroups of the general affine group.
In 2010, Bienert and Klopsch  \cite{BienertKlopsch}, studied the automorphism group of cyclic codes. They constructed two families of cyclic codes with prescribed automorphism groups and classified all cyclic codes with primitive automorphism groups. 
 In 2013, Ghorpade and Kaipa \cite{GhorpadeKaipa} determined the automorphism group of Grassmann codes, as well as related linear codes such as affine Grassmann codes and Schubert codes.
That same year, Guenda and Gulliver \cite{GuendaGulliver} investigated the Sylow $p$-subgroup of the automorphism group of cyclic codes of length $p^m$. 
  In 2021, Geiselhart et. al.  \cite{GEECB} introduced a new and effective type of permutation group decoder, known as the automorphism ensemble (AE) decoder, based on the automorphism group of given codes. 
Geiselhart et al. \cite{GEECB2} discussed the automorphism group of Polar codes, and proved that the block lower-triangular affine group  is contained in the automorphism group of polar codes. Subsequently,  it was proved by Li et al. \cite{Lietcl} that the block lower-triangular affine group is the full affine automorphism group of polar codes.

In this paper, we investigate the automorphism group of binary cyclic codes, and the paper is organized as follows. In Section 2, we briefly review the concepts of binary cyclic codes, permutations groups, code automorphisms and several useful known results. In Section 3, for cyclic codes of length $kn$ where $k, n$ are integers, a relatively large subgroup of automorphisms is characterized by a matrix presentation method of codewords.  
In Sections 4 and 5, respectively, we give classifications of full automorphism groups of  cyclic codes of lengths $2p$ and $p^n$ with irreducible generator polynomials, where $p$ is odd prime. 
Finally, we summarize our results and propose several questions for future consideration.

\section{Preliminaries}
Unless otherwise specified, all polynomials and matrices are defined over $\mathbb F_2$, the finite field with $2$ elements.

\subsection{Binary Cyclic Codes}

A binary {\em linear code} of length $n$ and dimension $k$, denoted by $[n, k]$-code, is a linear subspace $\mathcal C$ with dimension $k$ of the vector space $\mathbb{F}_2^n$. The vectors in $\mathcal C$ are called {\em codewords}. 
The {\em weight} of a codeword is the number of non-zero elements.

A binary linear code $\mathcal C$ of length $n$ is a cyclic code if every (right) cyclic shift of a codeword is also a codeword: 
 $$\forall c= (c_0, c_1, \cdots, c_{n-1})\in \mathcal C \Rightarrow   (c_{n-1}, c_0,  \cdots, c_{n-2})\in \mathcal C.$$
 As is well known that binary cyclic codes can be linked to ideals in polynomial rings. Let $R= \mathbb{F}_2[x]/(x^n-1)$ be a polynomial ring over $\mathbb{F}_2$. 
Define the polynomial associate to the codeword $c= (c_0, c_1, \cdots, c_{n-1})$ to be $c(x)= c_0+ c_1x+ \cdots +c_{n-1}x^{n-1}\in \mathbb{F}_2[x]$. Note that, $x\cdot c(x)= c_0x+c_1x^2+ \cdots +c_{n-1}x^{n}$, corresponds to a right cyclic shift of the polynomial with the representation $x\cdot c(x)= c_{n-1}+ c_0x+c_1x^2+ \cdots+ c_{n-2}x^{n-1}$. 
Then $\mathcal C$ is a principle ideal in $R$. 
The ideal is generated by the unique monic polynomial $g$ in $\mathcal C$ of minimum degree, which is called the {\em generator polynomial}. 
And thus for each polynomial $g(x)\in \mathbb{F}_2[x]$ divides the polynomial $x^n-1$ determines a cyclic $[n, n-k]$ code $\mathcal C$ with $k$ being the degree of $g(x)$.

\subsection{Permutation groups }

A  {\em permutation} $\pi$ on finite set $\Omega= \{1, 2, \cdots, N\}$ is a bijection from $\Omega$ to itself, and a finite permutation group $G$ on $\Omega$ is a finite group whose elements are permutations  with composition as the  group operation. 
Especially, the {\em symmetric group} of degree $N$, say ${\rm S}_{N}$,  is the set of all permutations $\pi$.
For each element $a\in \Omega$, the subgroup $G_a\, =\{ g\in G: g(a) = a\}$ of $G$ is called the {\em stabilizer} of $a$ and consists 
all permutations of $G$ fixed $a$, i.e. that permute $a$ to itself. 
A permutation group $G$ is called {\em transitive} if it acts transitively on set $\Omega$, i.e. for any pair of elements $x, y\in \Omega$ there exists a permutation $g\in G$ such that $g(x)= y$. And $G$ is called {\em 2-transitive} if for each element $a\in \Omega$, if the stabilizer $G_a$ is transitive on $\Omega\setminus \{a\} $. 
 The permutation group $G$  is said to be {\em regular} if $G$ is transitive and for all $a\in\Omega$, the stabilizer $G_a$ is trivial.  
For a transitive permutation group $G$ on $\Omega$, let $\mathcal P$ be a partition of $\Omega$. Then $\mathcal P$ is called {\em $G$-invariant} if for all $g\in G$ and $C\in \mathcal P$, the set $C^g= \{ g(a)\, |\, a\in C\} $ is belonging to $\mathcal P$. 
 A permutation group $G$ is called {\em primitive} on $\Omega$ if the only $G$-invariant partitions of $\Omega$ are the trivial ones, that is the universal partition $\{ \Omega\}$ and the partition consisting singleton subsets.

A finite permutation group which contains a cyclic regular subgroup is called a {\em c-group}.
Characterizing c-groups is an old topic in permutation group theory, initiated by Burnside (1900), and studied by Schur et al., see for example \cite[Chapter 4]{Wielandt}. 
A classical result of Schur tells that a primitive c-group is 2-transitive or has prime degree.  
Based on the finite simple group classification, a precise list of primitive c-groups was independently given in \cite{Jones, Li}.

\subsection{Code automorphisms}

For a code $\mathcal{C}$ of length $n$, a {\em code automorphism} $\tau$ of $\mathcal{C}$ is a permutation on the $n$ codeword bits that maps each codeword $x\in \mathcal{C}$ into another codeword $x'=\tau(x)\in \mathcal{C}$. The {\em automorphism group} of $\mathcal{C}$, denoted by $\Aut{\mathcal{C}}$,  is the set containing all automorphisms of code $\mathcal{C}$ using composition as the multiplication. 
Hence, if $\mathcal{C}$ is a cyclic code, one can see that the cyclic permutation $\tau= (1, 2, \cdots, n)$ is an automorphism, i.e. $\tau\in \Aut\mathcal{C}$. Moreover, $\tau$ generates a cyclic regular subgroup of $\Aut\mathcal{C}$, which implies that $\Aut\mathcal{C}$ is a c-group. 
Hence, in order to classify the full automorphism group of each cyclic code, one can firstly classify all c-groups. 

In \cite{BienertKlopsch}, base on a characterization of primitive c-group, Bienert and Klopsch gave the following characterization on binary cyclic codes with primitive automorphism groups. 

 \begin{theorem}\cite[Theorem E]{BienertKlopsch}{\label p} Let $G\leqslant {\rm S}_N$ be the automorphism group of a binary cyclic code ${\mathcal C}$, and suppose that $G$ is a primitive permutation group. Then one of the following holds. 
 
{\rm (1)} $C_p\lneqq G\lneqq {\rm AGL}(1, p) $ where $p=N\geq 5$ is a prime. 

{\rm (2)} $G= {\rm S}_N$; in this case $\mathcal C$ is one of four elementary codes.

{\rm (3)} $G= {\rm P\Gamma L}(d, q) $ where $d\geq 3, q=2^k$ for $k\in {\mathbb N}$ and $N= (q^d-1)/(q-1)$.

{\rm (4)} $G = M_{23}$ and $N= 23$. 

Moreover, each of the groups listed in {\rm (2)} to {\rm (4)} does occur as the automorphism group of a suitable binary cyclic code. 

 \end{theorem}

Note that, in Theorem~\ref{p}, the ${\rm S}_N$ denotes the symmetric group of degree $N$, the $C_p$ denotes the cyclic group of order $p$, the ${\rm AGL}(1, p)$ denotes the 1-dimensional affine general linear group of degree $p$, the ${\rm P\Gamma L}(d, q) $ denotes the projective semi-linear group, and
the $M_{23}$ denotes the Mathieu group of degree 23.

When the binary cyclic codes have a prime length $p$, the following corollary can be derived.    

    \begin{corollary} Let $G\leqslant {\rm S}_p$ be the full automorphism group of a binary cyclic code ${\mathcal C}$ of length $p$ prime. Then $G$ is primitive, and one of the following holds. 
 
{\rm (1)} $C_p\lneqq G\lneqq {\rm AGL}(1, p) $ where $p\geq 5$ is a prime; 

{\rm (2)} $G= {\rm S}_ p$; 

{\rm (3)} $G = M_{23}$ and $p= 23$. 
 \end{corollary}

  \section{Cyclic codes of length $kn$}

For positive integers $n$ and $k$, the polynomial $x^{kn}-1$ can be factorized into $ (x^n-1)(x^{(k-1)n}+x^{(k-2)n}+\cdots +x^{n}+1)$. Thus each factor of polynomial $x^n-1$ is also a factor of $x^{kn}-1$. 
Now, suppose that a polynomial $g(x)$ of degree $m>1$ is a factor of $x^n-1$. 
Denote $\mathcal{C}_n$ and $\mathcal{C}_{kn}$ respectively as binary cyclic codes of lengths $n$ and $kn$ defined both by $g(x)$ as the generator polynomial.
Let $v_0\in\mathcal{C}_n$ and $v\in\mathcal{C}_{kn}$ respectively be the generator codewords (row vector) corresponding to $g(x)$. One can see that the first $n$ columns of $v$ is actually $v_0$, and the last $kn-n$ columns are all zeros. 
Rewrite each codeword $c \in C_{kn}$ as a $k\times n$ matrix $M_c$ in a way where the $(in-n+1)$-th to the $in$-th codeword bits of $c$ being the $i$-th row of $M_c$ where $1\leqslant i\leqslant k$. Then the cyclic code $\mathcal{C}_{kn}$ can be viewed as a set of $k\times n$ matrices. For example, the codesword
$v$ of $\mathcal{C}_{kn}$ can be represented by the $k\times n$ generator matrix 
$$M_{v}= \begin{pmatrix}
v_0\\
o_1 \\ 
o_2\\
\cdots \\  
o_{k-1}  
\end{pmatrix}$$ 
where $o_i$ with $1\leqslant i\leqslant k-1$ represents a row vector of $n$-zeros. 

Then base on this matrix representation, one can establish the following results. 

\begin{lemma}\label{line} Each permutation on the $k$ rows results in a code automorphism of $\mathcal{C}_{kn}$, generating a subgroup that isomorphic to the symmetric group ${\rm S}_k$.
\end{lemma}
\pf 
Let 
 $c_1, c_2, \cdots, c_k\in {\rm S}_{kn}$ be $k$ cyclic permutations of order $n$ such that $c_1= (1, 2,\cdots, n)$, $c_2= (n+1, n+2,\cdots, 2n), \cdots, c_k=((k-1)n+1, (k-1)n+2,\cdots, kn) $.
 For $0\leqslant t\leqslant n-m-1$, let 
 $$M_{1,t}= \begin{pmatrix}
v_0^{c_1^t}\\
o_1 \\ 
o_2\\
\cdots \\  
o_{k-1}  
\end{pmatrix}, M_{2,t}=\begin{pmatrix}
o_1 \\ 
v_0^{c_2^t}\\
o_2\\
\cdots \\  
o_{k-1}  
\end{pmatrix}, \cdots,  M_{k,t}=\begin{pmatrix}
o_1 \\ 
o_2\\
\cdots \\  
o_{k-1}  \\
v_0^{c_k^t}
\end{pmatrix}$$ 
be generator matrices of ${\cal C}_{kn}$. 
 Then 
 the cyclic permutation 
$$\tau= (1, n+1, \cdots, (k-1)n+1)(2, n+2, \cdots, (k-1)n+2)\cdots (n, 2n, \cdots, kn)$$ permutes all these matrices by 
 acting on the $k$ rows.  
For any other generator matrices, one can observe that there are at most two rows containing 1s, while the rest rows are all zeros. Additionally, each column has a weight of at most one. 
Hence, $\tau$ is an automorphism of $\mathcal{C}_{kn}$. 

Moreover, if $k\geqslant 4$, by selecting a permutation that 
swaps two rows with all zeros, along with $\tau$, they generate a subgroup that is  isomorphic to the symmetric group ${\rm S}_k$. 
When $k= 2$, one can see that swapping the two rows is an automorphism as it is equivalent to the involution permutation $(1, n+1)(2, n+2)\cdots (n, 2n)$ which equals to the permutation $(1, 2,\cdots, 2n)^n$. 
When $k= 3$, let $$M= \begin{pmatrix}
v_1 \\ 
v_2\\
0
\end{pmatrix} \ {\rm and}\ \, M'= \begin{pmatrix}
v_1 \\ 
0\\
v_2
\end{pmatrix}$$
be two generator matrices of ${\cal C}_{kn}$ in which  
both $v_1$ and $v_2$ are non-zero. Note that, $M'$ is obtained by interchanging the second and third rows of $M$.  
Then one has $$M'= M+  M_1+M_2\,\  {\rm where}\ \, M_1= \begin{pmatrix}
0\\
v_1 \\ 
v_2
\end{pmatrix}\ {\rm and }\  \, M_2= \begin{pmatrix}
0\\
v_1+ 
v_2\\
0
\end{pmatrix}. $$ Note that if $M_2 \in \mathcal{C}_{kn}$, then $M'\in \mathcal{C}_{kn}$. 
In fact, for $v_1+v_2$, there exists an integer $\sigma$ such that $v_1+v_2= {v_0}^{{c_1}^\sigma}$, hence $M_2$ is a linear combination of 
$$\begin{pmatrix}
0\\
{v_0}^{c_2}\\
0
\end{pmatrix}, \begin{pmatrix}
0\\
{v_0}^{{c_2}^2}\\
0
\end{pmatrix}, \cdots,
\begin{pmatrix}
0\\
{v_0}^{{c_2}^{n-m}}\\
0
\end{pmatrix},
$$ 
and thus $M_2\in \mathcal{C}_{kn}$. Since $M_1\in \mathcal{C}_{kn}$, it follows that $M\in \mathcal{C}_{kn}$. Hence, along with $\tau$, ${\cal C}_{kn}$ admits an automorphism subgroup isomorphic to ${\rm S}_3$.  This finishes the proof of Lemma~\ref{line}.\qed

Furthermore, by considering the permutations on columns, one can obtain the following result. 
\begin{lemma}\label{column} The full automorphism group $\Aut\mathcal{C}_{kn}$ of $\mathcal{C}_{kn}$ contains a subgroup that is isomorphic to  ${\rm S}_k \times G$ where $G\cong \Aut\mathcal{C}_n $, the full automorphism group of $\mathcal {C}_n$.
\end{lemma}
\pf
For each permutation $\tau\in {\rm S}_n$, let $\bar\tau\in {\rm S}_{kn}$ be a lifting permutation of $\tau$  such that $\bar\tau= \tau{\rho_n}(\tau){\rho_{2n}}(\tau)\cdots {\rho_{(k-1)n}}(\tau)$ where $\rho_{in}$ maps coordinates $j$ to $in+j$ with $1\leqslant j\leqslant n$ and $1\leqslant i\leqslant k-1$. Then 
we claim that $\bar\tau\in\Aut{\mathcal C}_{kn}$ if and only if  $\tau\in\Aut{\mathcal C}_{n}$. 
First of all, it can be seen that if $\bar\tau\in\Aut{\mathcal C}_{kn}$ then $\tau\in\Aut{\mathcal C}_{n}$, we need to show the necessary condition. 

Now, assume that $\tau\in \Aut{\mathcal C}_{n}$. Similar to the arguments in Lemma~\ref{line}, one only needs to consider the action of $\bar\tau$ on generator matrices with two non-zero rows. Let 

$$M= \begin{pmatrix}
v_1\\
v_2\\
0\\
\cdots \\  
0  
\end{pmatrix}\ \,  {\rm and}\ \,  M^{\bar\tau}= \begin{pmatrix}
v_1^{\tau}\\
v_2^{\rho_n(\tau)}\\
0\\
\cdots \\  
0  
\end{pmatrix}$$
where  $v_1+v_2= v_0^{c^{t_0}}$ with $n-m-1<t_0<n$ and $c= (1, 2, \cdots, n)$. 
Then for $M^{\bar\tau}$, 
we have $v_1^{\tau}+v_2^{\rho_n(\tau)}= v_0^{c^{t_0}\cdot\tau}$. 
Note that, $v_0^{c^{t_0}\cdot\tau}\in \mathcal{C}_n$ as ${c^{t_0}\cdot\tau}\in \Aut\mathcal{C}_n$, and furthermore,  
it can be linearly represented by $v_0, v_0^{c},\cdots, v_0^{c^{n-m-1}}$. Hence, $M^{\bar\tau}\in \mathcal{C}_{kn}$. 
Similarly, by choosing different $M$, one can verify that $\bar\tau$ is an automorphism of $\mathcal {C}_{kn}$.  

Now one can establish that all the automorphisms $\bar{\tau}$  that permute columns generate a subgroup $G$ that is isomorphic to $\Aut\mathcal{C}_n$.  
In addition, this subgroup is commutative with respect to the actions of permuting rows and columns. Thus, along with Lemma~\ref{line}, one can deduce that $\Aut\mathcal{C}_{kn}$ contains a subgroup of automorphisms isomorphic to ${\rm S}_k\times G$. \qed

In the previous discussions, Lemma 3.1 and 3.2 have provided us with a fundamental comprehension of code automorphisms as complete row and column permutations. Naturally, we seek to extend these observations detailed consideration on row permutations and prove the following theorem.

\begin{theorem}\label{cl}  The full automorphism group $\Aut\mathcal{C}_{kn}$ of $\mathcal C_{kn}$ contains a subgroup that is isomorphic to  $({\rm S}_k)^n : G$ where $G\cong \Aut\mathcal{C}_n $, the full automorphism group of $\mathcal {C}_n$.

\end{theorem}
\pf By Lemma~\ref{line} and \ref{column}, one needs to show that each row permutation, say $\sigma$,  of a given non-zero column on generator matrices is an automorphism of $\mathcal C_{kn}$.  
Let $M$ be a generator matrix of $\mathcal C_{kn}$, and suppose $\sigma$ is the cyclic permutation $(b, n+b, \cdots, (k-1)n+b)$ of order $k$, which acts on a chosen non-zero $b$-th column. It can be observed  that the $b$-th column of matrix $M+M^{\sigma}$ has two 1s and all other columns are zeros. 
Next, we will show that the matrix 
$$A= \begin{pmatrix}
1& 0& \cdots& 0\\
1&  0&  \cdots& 0\\
0& 0&  \cdots & 0\\
&&      \ddots & \\  
0& 0&   \cdots&0\\
\end{pmatrix}$$ can be linearly represented by the generator matrices.

In fact, in $\mathcal{C}_n$, the $n-m$ number of generator vectors can be presented by $v_0, v_0^c, \cdots, v_0^{c^{n-m-1}}$ with $c= (1, 2, \cdots, n)\in \Aut\mathcal{C}_n$. Consequently, the vector $v_0^{c^{n-m}}$ can be linearly represented by $v_0, v_0^c,  \cdots$, $v_0^{c^{n-m-1}}$ as well. Thus there exists $k_1, k_2, \cdots, k_{n-m-1}\in \mathbb{F}_2$ such that $$v_0^{c^{n-m}}= v_0+ k_1v_0^c+k_2v_0^{c^2}+ \cdots + k_{n-m-1} v_0^{c^{n-m-1}}.$$ 
Note that $k_1, k_2, \cdots, k_{n-m-1}$ cannot be all zero as $v_0\ne v_0^{c^{n-m}}$. Let $\tilde c= (1, 2,\cdots, kn)\in \Aut\mathcal{C}_{kn}$ be a cyclic automorphism of order $kn$. 
Hence, we obtain the expression $$A={M_0}+k_1{M_0}^{\tilde c}+k_2{M_0}^{\tilde c^2}+\cdots+ k_{n-m-1}{M_0}^{{\tilde c}^{n-m-1}} +{M_0}^{{\tilde c}^{n-m}},$$ where

$$M_0= \begin{pmatrix}
v_0\\
0\\
\cdots \\  
0  
\end{pmatrix},$$
indicating that $A\in\mathcal C_{kn}$.  
Furthermore, 
there exists a proper integer $t_0$ such that $A^{{\tilde c}^{t_0}}= M+M^{\sigma}$, which implies that $M^{\sigma}\in\mathcal C_{kn}$ and thus 
$\sigma$ is a code automorphism of order $k$. Meanwhile, if $\sigma$ interchanges two zero entries in the $b$-th column of $M$, then it fixes $M$. 
Consequently, all such $\sigma$ generate a subgroup of automorphisms isomorphic to the symmetric group ${\rm S}_k$. Similarly, one can see that $\Aut\mathcal{C}_{kn}$ contains a subgroup isomorphic to $({\rm S}_k)^n$, which is generated by $S = \{ (i, n+i,  \cdots, (k-1)n+i), (i, n+i)\ |\ 1\leqslant i\leqslant n \}$. 
Moreover, by Lemma~\ref{line} and \ref{column}, one can see that the automorphism subgroup $G$ normalizes but not commute with $({\rm S}_k)^n$. Therefore, $\Aut\mathcal{C}_{kn}$ contains a subgroup of automorphisms isomorphic to $({\rm S}_k)^n : G$ with $G\cong  \Aut\mathcal{C}_n$, the full automorphism group of $\mathcal {C}_n$.
\qed
 
In the above discussion, Theorem~3.3 has shown that any row permutation of a fixing column is a code automorphism. The remaining permutations belongs to the set ${\rm S}_{kn}\setminus (({\rm S}_k)^{n}: G)$, that is oblique permutations on the generator matrices. 
Nonetheless, determining these permutations for arbitrary integers $n$ and $k$ may be challenging.  
In the following sections, we will determine the full automorphism group for particular integers  $n$ and $k$.

\section{Cyclic codes of length $2p$}

For a prime number $p>2$ and  the polynomial $x^{2p}-1$, we know that $x^{2p}-1= (x^p-1)^2$. This indicates that each irreducible factor of $x^{2p}-1$ is the same as $x^p-1$ but has a multiplicity of 2. 
Let ${\mathcal C}_{p, g}$ be a binary cyclic code of length $p$ that is generated by generator polynomial $g= g(x)$. The full automorphism group of ${\mathcal C}_{p, g}$, say $\Aut{\mathcal C}_{p, g}$, is of degree $p$, thus $\Aut{\mathcal C}_{p, g}$ is primitive and falls under the classification in Corollary~2.2. 

In this section, we investigate the full automorphism group of cyclic codes ${\mathcal C}_{2p, f}$ of length $2p$ that are generated by generator polynomial $f$ which equals to either $g$ or $g^2$ with $g= g(x)$ being irreducible factors of $x^{2p}-1$.  
We present the following results.

 \begin{lemma}\label{two} 
 Let ${\mathcal C}_{2p, g^2}$ be a binary cyclic code of length $2p$ with generator polynomial  $g^2$ where $g$ is irreducible. It follows that the full automorphism group $\Aut{\mathcal C}_{2p, g^2}$ of ${\mathcal C}_{2p, g^2}$ is isomorphic to $(G\times G). C_2$ where 
 $G\cong \Aut{\mathcal C}_{p, g}$. 
 \end{lemma}
 \pf Let $v$ be a codeword of the binary cyclic code ${\mathcal C}_{p, g}$ that corresponds to the polynomial $g$, and let $\tilde{v}$ be a codeword of ${\mathcal C}_{2p, g^2}$ that corresponds to the polynomial $g^2$. By rewriting $\tilde{v}$ in a format where all odd-indexed coordinates are listed in the first row, and all even-indexed coordinates are listed in  the second row, for example one can express 
$\tilde v$ as $\begin{pmatrix}v\\0\end{pmatrix}$. And in the following  for simplicity, we will say $\tilde v =  \begin{pmatrix}v\\0\end{pmatrix}$.

 For regular cyclic automorphisms $\tilde c= (1, 2, \cdots, 2p)\in \Aut{\mathcal C}_{2p, g^2}$ and $c= (1, 2, \cdots, p)\in \Aut{\mathcal C}_{p, g}$, it holds that 
 ${\tilde v}^{{\tilde c}}= \begin{pmatrix}0\\v\end{pmatrix}$ and ${\tilde v}^{{\tilde c}^2}= \begin{pmatrix}v^c\\0 \end{pmatrix}$. Furthermore, 
  ${\tilde v}^{{\tilde c}^t}$ equals to $ \begin{pmatrix} v^{c^{\frac{t}{2}}}\\ 0 \end{pmatrix}$ when $t$ is even,  or $ \begin{pmatrix}0\\v^{c^{\frac{t-1}{2}}}\end{pmatrix}$ when $t$ is odd. 
Additionally, for each codeword $w= \begin{pmatrix}w_1\\ w_2\end{pmatrix}\in \mathcal C_{2p, g^2}$, one can see that  $w_1, w_2\in \mathcal C_{p, g}$. 
  
 Hence, similar to the arguments in Lemma~\ref{column}, it follows that any permutation that  acts on a single line of $\tilde v$ is a code automorphism if and only if its restriction on odd (or respectively even) coordinates is a code automorphism of ${\mathcal C}_{p, g}$, which generates a subgroup isomorphic to $G\times G$ with $G\cong \Aut {\mathcal C}_{p, g}$. Furthermore, the permutation $\tau=(1, 2)(3,4)\cdots (2p-1, 2p)$ is clearly a code automorphism of ${\mathcal C}_{2p, g^2}$, which interchanges $G\times G$. 
Additionally, it can be observed that any other permutations not contained in the generating subgroup by $\tau$ and $G\times G$, which is isomorphic to $(G\times G).C_2$, cannot be  code automorphisms. 

 On the contrary, consider $\varphi\in{\rm S}_{2p}\backslash (G\times G).C_2$. The group $(G\times G).C_2$ acts transitively on the odd coordinates, we may assume that $\varphi$ fixes the first coordinates and flips certain other odd coordinates to even coordinates. Taking generator matrix $$M=  \begin{pmatrix}
v\\
0 
\end{pmatrix}= \begin{pmatrix}
1& v_{1}\\
0& 0 
\end{pmatrix}, {\rm and\   suppose\ that\   }
M^{\varphi}= \begin{pmatrix}
1&v_1'\\
*&*
\end{pmatrix}\ne M $$ where $v_1'$ is obtained from $v_1$ by exchange some 1s to 0s. 
Hence, $M^{\varphi}\notin\mathcal C_{2p, g^2}$, otherwise $(1, v_1')$ has to be contained in $\mathcal C_{p,g}$, a contradiction. Thus $\varphi$ cannot be an automorphism. 
  Therefore, the Lemma holds. \qed

 Now one can get the characterization of full automorphism group of binary cyclic codes of length $2p$ as follows. 

\begin{theorem}\label{2p} The full automorphism group $\Aut{\mathcal C}_{2p, f}$ is isomorphic to one of the following: 

{\rm (1)} ${\rm S}_{2p}$, with $f= x+1$; or 
 
{\rm (2)} $({\rm S}_2)^p: \Aut{\mathcal C}_{p, f}$ with $f$ irreducible and not equal to $x+1$; or  

{\rm (3)} $(\Aut{\mathcal C}_{p, g}\times \Aut{\mathcal C}_{p, g}). C_2$ with $f= g^2$ and $g= g(x)$ being irreducible. 
\end{theorem}
\pf According to Lemma~\ref{line}, ~\ref{cl} and ~\ref{two}, we only need to consider case (2) that beyond $({\rm S}_2)^p: \Aut{\mathcal C}_{p, g}$ there exists no additional code automorphisms. 
On the contrary, let $\varphi\in{\rm S}_{2p}\backslash (({\rm S}_2)^p: \Aut\mathcal C_{p, g})$. As the group $({\rm S}_2)^p: \Aut{\mathcal C}_{p, g}$ acts transitively on $\Omega= \{1, 2, \cdots, 2p\}$, we may assume that $\varphi$ fixes the $(1, 1)$-entry and maps the $(2, 1)$-entry of each generator matrix to the $(1, 2)$-entry. By arguments similar to those in  Lemma~\ref{cl}, we know that there exists a codeword 

$$A= \begin{pmatrix}
1& 0& \cdots& 0\\
1&  0&  \cdots& 0 
\end{pmatrix}. $$
In particular, the image of $A$ under $\varphi$ is 

$$A^{\varphi}= \begin{pmatrix}
1& 0& \cdots& 0\\
1&  0&  \cdots& 0 
\end{pmatrix}^{\varphi} =  \begin{pmatrix}
1& 1&0& \cdots& 0\\
0&  0&0&  \cdots& 0 
\end{pmatrix}. $$
Note that $A^{\varphi}\notin \mathcal C_{2p, f}$ when $f\ne x+1$. Therefore  $\varphi\notin \Aut\mathcal C_{2p, f}$, and this finishes the theorem.  
\qed

\section{Cyclic codes of length $p^n$}
For prime $p\ne 2$, the polynomial $x^{p}-1$ satisfies that $x^{p}-1= 
(x+1)
(x^{p-1}+x^{p-2}+\cdots+ x+1)$. Let $f(x)= x^{p-1}+x^{p-2}+\cdots+ x+1$. Then we have 
\begin{align*}
x^{p^n}-1&= \ \ 
(x^{p^{n-1}}+1)f(x^{p^{n-1}})\\
&= \ \ 
(x^{p^{n-2}}+1)f(x^{p^{n-2}})f(x^{p^{n-1}})\\
&= \ \ \ \  \cdots \\
&= \ \ (x+1)f(x)f(x^p)\cdots f(x^{p^{n-2}})f(x^{p^{n-1}}). 
\end{align*}
%
Now, suppose that $f(x)= f_1(x)f_2(x)\cdots f_s(x)$ is the irreducible factorization with $f_i(x)$ being irreducible. Then, $f_i(x^p),f_i(x^{p^2}), \cdots, f_i(x^{p^{n-1}})$ are all irreducible.

Now let $\mathcal C_{p^n,\, m}$ be a binary cyclic code of length $p^n$ with generator polynomial $g= f_i(x^{p^m})$ for $0\leqslant m\leqslant n-1$. Let $v_{m,\, n}$ be the codeword corresponding to $g$. 
 In particular, when $m= 0$, let $v_{0,\, n}\in \mathcal C_{p^n,\, g}$ be the codeword of length $p^n$ corresponding to the generator polynomial $f_i(x)$. 
We represent $v_{m,\, n}$ by a $(p^{m}\times p^{n-m})$ matrix $M_{m,\, n}$ such that the  $(a, b)$-entry of $M_{m,\, n}$ is the  $(a+(b-1)p^{m})$-th coordinate of $v_{m,\, n}$. 
For example, when $m= 0$, we have 
$M_{0,\, n}=  v_{0,\, n}$, and 
$$M_{1,\, n}= \begin{pmatrix}
v_{0,\, n-1}\\
0\\
\cdots\\
0
\end{pmatrix}_{p\times p^{n-1}}\ {\rm when}\ m=1,\ {\rm and\ we\ have}\ M_{2,\, n}= \begin{pmatrix}
v_{0,\, n-2}\\
0\\
\cdots\\
0
\end{pmatrix}_{p^2\times p^{n-2}}\ {\rm when}\ m= 2.$$  
In paricular, $$M_{n-1,\, n}= \begin{pmatrix}
v_{0,\, 1}\\
0\\
\cdots\\
0
\end{pmatrix}_{p^{n-1}\times p}\ {\rm  when}\ m= n-1.$$


Then we have the following results. 

 \begin{lemma} 
The full automorphism group $\Aut\mathcal C_{p^{n}, 0}$ of $\mathcal C_{p^{n}, 0}$ is isomorphic to $({\rm S}_{p^{n-1}})^p: \Aut\mathcal C_{p, 0}$ where $\Aut\mathcal C_{p, 0}$ is the full automorphism group of $\mathcal C_{p, 0}$. 
\end{lemma}
\pf
According to Theorem~\ref{cl}, it suffices to show that besides $({\rm S}_{p^{n-1}})^p: \Aut\mathcal C_{p, 0}$ there are no additional automorphisms. Similar to the arguments in Lemma~\ref{2p}, one can verify that the lemma holds. 
\qed 

 \begin{theorem} 
For $m\geqslant 1$, the full automorphism group $\Aut\mathcal C_{p^n,\, m}$ of $\mathcal C_{p^{n}, m}$ is isomorphic to $(\Aut\mathcal C_{p^{n-m}, 0})^{p^{m}}: {\rm S}_{p^{m}}$ where $\Aut\mathcal C_{p^{n-m}, 0}$ is the full automorphism group of $\mathcal C_{p^{n-m}, 0}$. 

\end{theorem}
\pf 
Let $M_{m, n} \in \mathcal C_{p^n,\, m}$ be a $p^m\times p^{n-m}$ matrix such that $$ M_{m, n}=\begin{pmatrix}
v_{0,\, n-m}\\
0\\
\cdots\\
0
\end{pmatrix}.$$ Then for an automorphism ${\bar c}= (1, 2, \cdots, p^n)\in \Aut\mathcal C_{p^n,\, m}$, 
 one has that the image of ${M_{m, n}}$ under the permutation action of ${\bar c}$ gives 
$${(M_{m, n})}^{\bar c}= \begin{pmatrix}
0\\
v_{0,\, n-m}\\
0\\
\cdots\\
0
\end{pmatrix}\in \mathcal C_{p^n,\, m}.$$ In fact, each generator matrix can be represented by the images of ${M_{m, n}}$ under the permutation action of ${\bar c}^{\kappa}$, which is  
$$(M_{m, n})^{{\bar c}^{\kappa}}= \begin{pmatrix}
0\\
\cdots\\
{v_{0,\, n-m}}^{c^{k}}\\
\cdots\\
0
\end{pmatrix}$$ where $c= (1, 2, \cdots, p^{n-m})\in\Aut\mathcal C_{p^{n-m},\, 0}$ and $\kappa= kp^m+t$ with $0\leqslant t\leqslant p^m$, in particular, ${v_{0,\, n-m}}^{c^{k}}$ lies in the $t$-th row. 
Hence, for each codeword $$M= \begin{pmatrix}v_0\\ v_1\\
\cdots\\
v_{p^m-1}\end{pmatrix}\in \mathcal C_{p^n, m},$$ one can see that  $v_{j}\in \mathcal C_{p^{n-m}, 0}$ with $0\leqslant j\leqslant p^m-1$. 
 Then similar to the arguments of Lemma~\ref{two}, one can verify this Theorem. 
\qed



\section{Conclutions}
This paper presents a matrix technique for describing long-length cyclic codes from given ones, enabling efficient identification of automorphism subgroups. Specifically, for any odd prime $p$, we utilize this method to determine the full automorphism groups of binary cyclic codes of lengths $2p$ and $p^n$ that are generated by irreducible generator polynomials. 
The classification of automorphisms for other cyclic codes remains open. Nevertheless, we believe that by a similar technique one can determine the automorphisms of cyclic codes of length $2p^n$ with irreducible generator polynomials as well. Consequently, we pose this as a open problem for future consideration. 

\begin{question} Classify full automorphisms of binary cyclic codes of length $2p^n$ with irreducible generator polynomial. 
\end{question}
 
Additionally, in our work, we did not consider much on generator polynomials which are not irreducible. With the help of {\sc Magma}, we have the following computational results in Table~1, and we would propose the second question for future consideration as well. 
 
 \begin{question} Classify automorphisms of binary cyclic codes of length $2p, p^n, 2p^n$ with $n\geqslant 2$. 
\end{question}

At last, one may consider binary cyclic codes of length $pq$ with $p$ and $q$ are both odd primes. The main reason for not considering this here is that the polynomial $x^{pq}-1$ has irreducible factors that are factors of neither $x^p-1$ nor $x^q-1$. 
For example, when $pq= 15$, $$x^{15}-1= (x^5+1)(x^2+x+1)(x^4+x+1)(x^4 + x^3 + 1)$$ where the later two factors are not factors of $x^3-1$. 
Hence, the classification of those cyclic codes may need the classification of permutation c-group of degree a product of two odd primes.

  \begin{center}{\small 
\captionof{table}{Binary cyclic codes $\mathcal C_n$ of length $n$}
\begin{tabular}{ cl c l  cl cl }
 \hline
$\bf n$ & Generator polynomial &  & $|\Aut\mathcal C_n|$ & \\ \hline
  $7$ & $x^3 + x + 1$  &  &168 &   \\ \hline     $7$ & $(x^3 + x + 1)(x^3 + x^2 + 1)$  && 7!& & \\ \hline
  $14$ & $(x^3 + x + 1)^2$ && $2\cdot 168^2$ \\ \hline   
    $14$ & $(x^3 + x + 1)(x^3 + x^2 + 1)$ & &$(7!)\cdot2^7$ \\ \hline   
        $49$ & $(x^3 + x + 1)(x^3 + x^2 + 1)$ & &$ (7!)^8$ \\ \hline  
         $49$ & $(x^{21} + x^7 + 1)(x^{21} + x^{14} + 1)$ & &$(7!)^8$ \\ \hline     
  $98$ & $(x^3 + x + 1)^2$ && $2\cdot 168^2\cdot (7!)^{14}$ \\ \hline   
 $98$ & $(x^3 + x + 1)(x^3 + x^2 + 1)$ & &$(7!)\cdot (14!)^7$ \\ \hline  

   $31$ & $x^5 +x^2+ 1$,  $x^5 + x^3 + 1$ or $x^5 + x^3 + x^2 + x + 1$ &&9999360 & & \\ \hline    
     $31$ & $(x^5 +x^2+ 1)(x^5 + x^3 + 1)$  & &310 & &  \\ \hline  
$31$ & $(x^5 +x^2+ 1)(x^5 + x^3 + 1)(x^5 + x^3 + x^2 + x + 1)$  & &155 & &  \\ \hline  
      $62$ & $(x^5 + x^2 + 1)^2$ & &$  2\cdot 9999360^2$ \\ \hline   $62$ & $(x^5 +x^2+ 1)(x^5 + x^3 + 1)$ & &$  310\cdot 2^{31}$ \\ \hline   
$62$ & $(x^5 +x^2+ 1)(x^5 + x^3 + 1)(x^5 + x^3 + x^2 + x + 1)$  & &$ 155\cdot 2^{31}$ & &  \\ \hline  
       $31^2$ & $(x^5 +x^2+ 1)(x^5 + x^3 + 1)$  & &$310\cdot (31!)^ {31}$& &  \\ \hline 
       $2\cdot31^2$ & $(x^5 +x^2+ 1)(x^5 + x^3 + 1)$  & &$ 310\cdot (62!)^{31}$& &  \\ \hline 

\end{tabular}}
 \end{center}

 \bigskip\bigskip
\noindent {\Large\bf Acknowledgments}

\smallskip
\noindent
The authors acknowledge the use of the {\sc Magma} computational package \cite{Magma}, which helped show the way to many of the results given in this paper.

\end{document}